%% file: main.tex
\documentclass{svproc}
\usepackage{amssymb}
\usepackage{graphicx}
\usepackage{tikz}
\usepackage{amsmath}
\usepackage{todonotes}
\usepackage{amsfonts}
\usepackage{url}
\usepackage{color}
\usepackage{xcolor}
\usepackage{array}
\usepackage{mathtools}
\usepackage{enumitem}
\usepackage{hyperref}
\usepackage{microtype}

\usepackage[sectionbib]{natbib}

\newcommand{\eg}{\textit{e.g.}}
\newcommand{\ie}{\textit{i.e.}}

\newcommand{\cbm}{\text{\,m}^3}
\newcommand{\kgpercbm}{\text{\,kg/m}^3}
\newcommand{\Land}{^{\text{land}}}
\newcommand{\Sea}{^{\text{sea}}}
\newcommand{\cont}{\text{cont}}
\newcommand{\FCL}{\,\text{FCL}}
\newcommand{\NVOCC}{\,\text{NVOCC}}

\newcommand{\CC}{C\nolinebreak\hspace{-.05em}\raisebox{.4ex}{\tiny\bf +}\nolinebreak\hspace{-.10em}\raisebox{.4ex}{\tiny\bf +}}
\newcommand{\Rge}{\mathbb{R}_{\ge 0}}
\DeclareMathOperator{\st}{\,s.t.}

\title{Partitioned vs.\ Integrated Planning of Hinterland Networks for LCL Transportation}
\titlerunning{LCL Hinterland Network Planning}
\author{Niklas Jost\inst{1} \and Dorothee Henke\inst{2} \and Ivo Hedtke\inst{3} \and Oliver Bredtmann\inst{3} \and Joachim Weise\inst{3} \and Christoph Buchheim\inst{2} \and Uwe Clausen\inst{1}}
\authorrunning{N.~Jost, D.~Henke, I.~Hedtke, O.~Bredtmann et al.}
\institute{%
Institute of Transport Logistics, TU Dortmund University, Germany, \email{niklas.jost@tu-dortmund.de}
\and%
Department of Mathematics, TU Dortmund University, Germany
\and%
Global Data Strategy \& Analytics, Schenker AG, Germany
}

\allowdisplaybreaks
\begin{document}
\bibliographystyle{agsm}

\maketitle

\begin{abstract}
Utilizing existing transportation networks better and designing (parts of) networks involves routing decisions to minimize transportation costs and maximize consolidation effects.
We study the concrete example of hinterland networks for the truck-transportation of less-than-container-load (LCL) ocean freight shipments:
A set of LCL shipments is given. They have to be routed through the hinterland network to be transported to an origin port and finally to the destination port via ship. On their way, they can be consolidated in hubs to full-container-load (FCL) shipments. The overall transportation cost depends on the selection of the origin port and the routing and consolidation in the hinterland network.
A problem of this type appears for the global logistics provider DB Schenker. We translate the business problem into a hub location problem, describe it mathematically, and discuss solution strategies. As a result, an integrated modeling approach has several advantages over solving a simplified version of the problem, although it requires more computational effort.
\keywords{hub location, logistics, algorithms, networks, discrete optimization, transport}
\end{abstract}

\section{Introduction}

Any international logistics provider has the task of building a transportation network involving different modes of operations. Especially for shipments to other continents, it can be necessary to use ships. In the transportation network, it needs to be decided which port is used for which shipment. Furthermore, in land transport, it is useful to consolidate multiple shipments. For instance, the task is to choose a hub (consolidation depot) and a port for any shipment, starting in Europe and ending on a different continent. An overview of this topic is given by \cite{notteboom2008bundling}. DB Schenker is an international logistics provider facing the challenge of building such a transportation network.

Routing shipments from an origin to a destination is similar to the minimum cost multi-commodity flow problem, where each of a number of different commodities must be transported from a source to a sink while satisfying flow constraints. In contrast to the problem at hand, the number of stopovers for a flow is not limited, a maximum capacity for each connection exists, and the cost for using a connection is linear in the volume. In the problem at hand, there is no limited capacity. In this case, the multi-commodity flow problem can easily be solved by any shortest path algorithm since each commodity flow is independent. Algorithms for the minimum cost multi-commodity flow problem can be found in \cite{awerbuch1994improved} or \cite{barnhart1994column}.

Another type of problems facing combined location and routing problems are hub location problems (HLP). All such problems have in common that each of several commodities needs to be routed between two branches using one or two hubs as stopovers. Using a hub-to-hub connection takes consolidation effects into account, which discounts the cost of this connection. Furthermore, two types of HLPs can be distinguished: single-allocation problems and multi-allocation problems. In single-allocation HLPs, every branch is connected to exactly one hub. Every shipment from or to this branch must be sent via the connected hub. In multi-allocation HLPs, hubs are chosen for every shipment separately. Further HLP variants differ regarding some constraints or objectives. A lot of research was done in this area including a branch-and-price approach by \cite{contreras2011branch}, a genetic algorithm by \cite{topcuoglu2005solving}, or a decomposition approach by \cite{kammerling2021dekompositionsverfahren}. A detailed list of different hub location problems and reviews can be found in \cite{campbell2012twenty,alumur2008network}. HLPs and the logistics provider's problem at hand have in common that every shipment needs to have one or two stopovers. In addition, in both problems the cost of a connection may be discounted by consolidation effects. Hence, the problem at hand can be understood to be a variant of the HLP. It also extends the classical HLP variants by a more realistic cost structure. 

The costs on every transport connection will be considered to be a stepwise function depending on the volume. This structure is mainly motivated by fixed costs for each container. A hub location problem with stepwise costs is investigated by \cite{RostamiBuchheim}. Based on this model, an approach for modeling and solving DB Schenker's hinterland HLP was constructed.

This paper considers the transportation of shipments not filling a truck (LTL, less-than-truck-load) and not filling an ocean freight container (LCL, less-than-container-load).
They can be consolidated in the land transport leg to FTL (full-truck-load) shipments.
They have to be consolidated into containers at the port using a CFS (container freight station) or NVOCC (non-vessel operating common carrier).
We use the terms LCL and NVOCC synonymously, and also differentiate between LTL and LCL only when necessary.

The paper is structured as follows. In Section~\ref{chap_challenges}, the challenges in the routing decision will be presented
from the point of view of the logistics provider DB Schenker. In Section~\ref{chap_math_model}, the concrete problem statement will be given, together with an explanation what is optimized. In Section~\ref{chap_schenker_sol}, the solution approach
of DB Schenker to this problem will be shown. This is a two-stage algorithm partitioning the problem
to make it tractable. A heuristic is used to solve the sub-problems of the stages. Hence, this method
had some downsides. Together with the Chair of Discrete Optimization and the Institute of Transport
Logistics of TU Dortmund University, the problem was faced again in a one-stage, integrated (non-partitioned)
manner. The resulting model is presented in Section~\ref{chap_tud_sol}. In Section~\ref{chap_results}, the
results are discussed, a conclusion to the methods is drawn and an outlook to further
research is given.
\input{chapters/1Challenges}
\input{chapters/2MathModel}
\input{chapters/3SchenkerSolution}
\input{chapters/4Hinterland}
\input{chapters/5Results}

\section*{Acknowledgments}
This work has been supported by Deutsche Forschungsgemeinschaft (DFG) under grants~BU 2313/6 and~CL 318/32-1.
We thank Carlos Rodriguez for doing part of the development and implementation of the model.

\bibliography{bib}

\end{document}

%% file: chapters/1Challenges.tex
\section{The hub location problem of a logistics provider%
\label{chap_challenges}}
In this section, the challenges in the routing decision will be presented from the point of view of the logistics provider DB~Schenker.
The routing of individual shipments through the transportation network is translated into a hub location problem.
This does not only answer how the shipments should be routed but also which entities of the existing network should be upgraded (introduction of hubs) in order to keep up with the volume and reduce cost by consolidation.
The logistics provider also provides the real-world instances from its Europe-APAC business, which will be discussed in Section~\ref{chap_schenker_sol}. The instance consists of around $400$ branches and hubs, more than $50$ origin ports, more than $300$ destination ports and several hundred thousands of shipments.

Transporting LCL shipments through NVOCCs in DB Schenker’s global Ocean Freight network is a complex challenge with various optimization opportunities.
The profitability of the business depends mainly on the utilization of its own containers and on the cost of land pre- and on-carriage.
Both factors are impacted by which origin port is chosen to consolidate cargo for a certain destination into full containers.
In addition to consolidating cargo at the origin port directly, there is also an option to introduce hinterland hubs for consolidation.

The problem presented in this paper, with a pan-European scope, is driven by DB Schenker's LCL Europe division with the intention to optimize the cost by a better usage of the European outbound LCL network.
The LCL organization establishes its own container services and maintains its own CFSs where shipments are loaded into ocean freight containers.
An integrated network optimization is presented that considers both options.
The main goal of this study can be summarized as follows:
\emph{minimizing the overall transportation cost by assigning catchment areas to origin ports (per destination) and by determining the optimal number and location of hinterland hubs.}

In more detail, the problem can be separated into two mutually dependent parts:
\begin{itemize}
\item \emph{Port assignment}: Find fixed assignments of catchment areas to origin ports that minimize overall transportation costs per destination.
\item \emph{Hinterland hubs}: Find optimal number and location of hinterland hubs.
\end{itemize}

\paragraph{Port assignment}
When there is enough load on an ocean relation, DB Schenker can operate own container services.
In this case, DB Schenker operates a CFS at the origin port and a CFS at the destination port where containers are unloaded.
Instead, when the load is low, it is given to consolidators (NVOCCs), which has the advantage of not requiring any own infrastructure or personnel at every port.
Typically, well-utilized own container services are much more cost effective.

In network design, this leads to a tradeoff between ocean and land transport costs.
Bringing cargo to very few big ports and consolidating it into own containers can bring down the ocean transport costs,
whereas using an NVOCC at a potentially more closely located port may lead to considerably lower land transport costs.

\paragraph{Hinterland hubs}

In order to reduce transportation cost in the pre-carriage (inland movement that takes place prior to the shipment being delivered to the port) stage, the consolidation of shipments within the hinterland network is considered.
Existing entities (\eg, origin branches) within the network can be upgraded to \emph{hubs} where shipments are consolidated and afterwards transported to their origin port.
Depending on the port assignment, a fixed assignment of hubs to origin branches has to be constructed.
When a new shipment enters the system at its origin branch with a specific destination port, the port assignment specifies which origin port to use and then there is only one decision left, namely whether it should be transported directly from origin branch to origin port, or to a hub where it gets consolidated with other shipments for that origin port.

\begin{figure}[b]
\newlength{\xunit}
\setlength{\xunit}{35mm}
\newlength{\yunit}
\setlength{\yunit}{7mm}
\newlength{\portnodewidth}
\setlength{\portnodewidth}{16mm}
\newlength{\terminalnodewidth}
\setlength{\terminalnodewidth}{20mm}
\centering
\begin{tikzpicture}
\useasboundingbox (\xunit,-1.75\yunit) -- (4\xunit,2\yunit);
\path [fill=black!15] (3\xunit - 0.5\portnodewidth,1.5\yunit) rectangle (3\xunit + 0.5\portnodewidth,-1.5\yunit);
\node at (3\xunit,1.5\yunit) [above] {origin port};
\node[align=center] (originbranch) at (1\xunit,0) {origin\\branch};
\node (terminal1) at (2\xunit,0.5\yunit) [minimum width=\terminalnodewidth] {hub~$H_1$};
\node (terminal2) at (2\xunit,-0.5\yunit) [minimum width=\terminalnodewidth] {\vdots};
\node (port1) at (3\xunit,1\yunit) [minimum width=\portnodewidth] {port~$P_1$};
\node (port2) at (3\xunit,0) [minimum width=\portnodewidth] {port~$P_2$};
\node (port3) at (3\xunit,-1\yunit) [minimum width=\portnodewidth] {\vdots};
\node[align=center] (destination) at (4\xunit,0) {destination\\port};
\draw[->] (port1.east) to[bend left] node[above,sloped] {\scriptsize CFS} (destination.north west);
\draw[->] (port2.east) to[bend left=20] node[above,sloped] {\scriptsize CFS} (destination.west);
\draw[->] (port2.east) to[bend right=20] node[below,sloped] {\scriptsize NVOCC} (destination.west);
\draw[->] (port3.east) to[bend right] node[below,sloped] {\scriptsize \ldots} (destination.south west);
\draw[dashed,->] (originbranch) -- (terminal1);
\draw[dashed,->] (originbranch) -- (terminal2);
\draw[dotted,->] (originbranch) to [bend left] (port1.west);
\draw[dotted,->] (originbranch) to (port2.west);
\draw[dotted,->] (originbranch) to [bend right] (port3.west);
\draw[dashed,->] (terminal1.east) -- (port1.west);
\draw[dashed,->] (terminal1.east) -- (port2.west);
\draw[dashed,->] (terminal1.east) -- (port3.west);
\draw[dashed,->] (terminal2.east) -- (port1.west);
\draw[dashed,->] (terminal2.east) -- (port2.west);
\draw[dashed,->] (terminal2.east) -- (port3.west);
\end{tikzpicture}
\caption{%
A shipment is transported from the origin branch to a destination port. An origin port is chosen and the shipment can be sent to the origin port directly (dotted lines) or via a hinterland hub (dashed lines). Origin port~$P_1$ has a CFS; at port~$P_2$ a CFS and an NVOCC can be used.
\label{fig:shipment_flow}}
\end{figure}
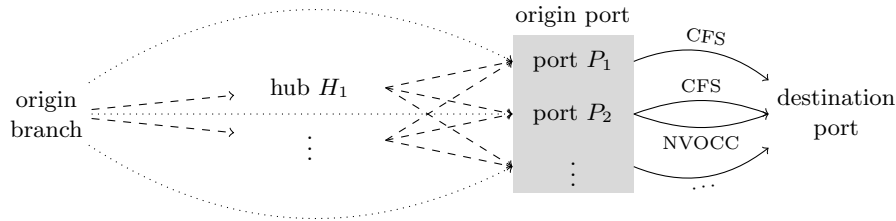

\paragraph{Shipment flow}
Figure~\ref{fig:shipment_flow} shows all possibilities for the flow of a single shipment through the hinterland and ocean network.
Shipments are always picked up at the consignor and transported to the nearest land transport branch within the hinterland network.
We can therefore assume that shipments originate directly at the hinterland branches as the consignor-branch leg does not depend on the presented network design process.
Furthermore, the transportation from the destination port to the consignee is not in the scope.
We therefore assume that shipments destinate at the destination port.
A shipment is transported either directly or via one hub to its origin port.
The origin port is not fixed as the decision which origin port should be used to transport the shipment to its fixed destination port is part of the problem.

The shipment flow presented above takes only a single shipment into account.
When multiple shipments are routed through the network simultaneously, an important requirement has to be respected:
Per pair $(\mathit{originbranch}, \mathit{originport})$, at most one hub can be used for the shipment flow, independently of the, possibly different, destination ports the shipments have.
Although it could reduce costs even further to allow multiple hubs, as in $\mathit{originbranch}\to \mathit{hub}_1 \to \mathit{originport} \to \mathit{destinationport}_1$ and $\mathit{originbranch}\to \mathit{hub}_2 \to \mathit{originport} \to \mathit{destinationport}_2$, it is practically not manageable to do so.

\paragraph{Cost structure}
Ocean transport costs are given per \emph{full container}.
NVOCC costs are given per \emph{volume}.
Land transport costs are given as cost per \emph{distance} per \emph{weight} and are stepwise constant for certain distance and weight intervals.

%% file: chapters/2MathModel.tex
\section{Mathematical Model} \label{chap_math_model}

The hub location problem
described in Section~\ref{chap_challenges}
will now be formulated more precisely and modeled mathematically.
Both the \emph{port assignment} and the \emph{hinterland hubs} decisions
are combined in an integrated modeling approach.

For a precise formulation of the problem,
the following data has to be given as the input:
\begin{itemize}
\item a finite set $B$ of branches,
\item a finite set $S$ of origin ports,
\item a finite set $T$ of destination ports,
\item demand volumes $v_{bt} \in \Rge$ for all $b \in B$ and $t \in T$,
\item cost functions $\Gamma\Land_{br} \colon \Rge \to \Rge$ for the hinterland transport, for all $b \in B$ and $r \in B \cup S$, depending on the transport volume,
\item cost functions $\Gamma\Sea_{st} \colon \Rge \to \Rge$ for the sea transport, for all $s \in S$ and $t \in T$, depending on the transport volume,
\item set-up cost $e_b$ for all $b \in B$,
\item consolidation cost $f_b$ per m$^3$ for all $b \in B$,
\item consolidation cost $g_s$ per m$^3$ for all $s \in S$.
\end{itemize}
More details on the cost functions $\Gamma\Land_{br}$ and $\Gamma\Sea_{st}$
will be discussed in Sections~\ref{chap_schenker_sol} and~\ref{chap_tud_sol}.
For now, we just consider them as general functions mapping the volume that is transported on the given transport connection to some cost.

In order to solve the problem in a feasible way,
the following decisions have to be made:
\begin{itemize}
\item For each transport connection from branch $b \in B$ to destination port $t \in T$: Select an origin port $s \in S$.
\item Select a set $H \subseteq B$ of branches as hubs.
\item For each (potential) transport connection from branch $b \in B$ to origin port $s \in S$:
Decide which percentage of the volume between $b$ and $s$ is sent via direct transport, and decide over which hub $h \in H$ the remaining volume is sent (and consolidated with other shipments there). If $b \in H$, only direct transport can be used.
\end{itemize}

When making these decisions,
the objective is to minimize the total cost,
which is the sum of:
\begin{itemize}
\item set-up cost for each hub $h \in H$,
\item consolidation cost for the flow routed via each hub $h \in H$,
\item consolidation cost for the flow routed via each origin port $s \in S$,
\item transportation cost for the direct hinterland transport between each branch $b \in B$ and origin port $s \in S$, together with the hub-based hinterland transport from $b$ to $s$ if $b \in H$,
\item transportation cost for the hub-based hinterland transport between each branch $b \in B$ and hub $h \in H$,
\item transportation cost for the sea transport between each origin port $s \in S$ and destination port $t \in T$.
\end{itemize}

We can also formulate the problem as a mathematical program using the following decision variables:

\begin{itemize}
\item $\boldsymbol z_{bts} \in \{0, 1\}$ for all $b \in B$, $t \in T$, and $s \in S$, determining if $s$ is used as origin port (1) or not (0) on the connection from $b$ to $t$,
\item $\boldsymbol x_b \in \{0, 1\}$ for all $b \in B$, determining if $b$ is a hub ($1$) or not ($0$),
\item $\boldsymbol y_{bs0} \in [0, 1]$ for all $b \in B$ and $s \in S$, determining which percentage of the volume on the connection from $b$ to $s$ uses direct transport,
\item $\boldsymbol y_{bsh} \in \{0, 1\}$ for all $b \in B$, $s \in S$, and $h \in B$, determining if $h$ is used as a hub ($1$) or not ($0$) on the connection from $b$ to $s$.
\end{itemize}

In terms of these variables,
the objective and constraints can be written as follows:
\begin{align}
\min \quad & \sum_{h \in B} e_h \boldsymbol x_h \label{1b_set-up} \\
& + \sum_{h \in B} f_{h} \sum_{b \in B} \sum_{s \in S} \sum_{t \in T} v_{bt} \boldsymbol z_{bts} \boldsymbol y_{bsh} (1 - \boldsymbol y_{bs0}) \label{1b_consolidation_hubs} \\
& + \sum_{s \in S} g_s \sum_{b \in B} \sum_{t \in T} v_{bt} \boldsymbol z_{bts} \label{1b_consolidation_ports} \\
& + \sum_{b \in B} \sum_{s \in S} \Gamma\Land_{bs}\!\left(\sum_{t \in T} \left(v_{bt} \boldsymbol z_{bts} \boldsymbol y_{bs0} + \sum_{c \in B} v_{ct} \boldsymbol z_{cts} \boldsymbol y_{csb} (1 - \boldsymbol y_{cs0}) \right)\right) \label{1b_transportation_branch_port} \\
& + \sum_{b \in B} \sum_{h \in B} \Gamma\Land_{bh}\!\left(\sum_{s \in S} \sum_{t \in T} v_{bt} \boldsymbol z_{bts} \boldsymbol y_{bsh} (1 - \boldsymbol y_{bs0}) \right) \label{1b_transportation_branch_hub} \\
& + \sum_{s \in S} \sum_{t \in T} \Gamma\Sea_{st}\!\left(\sum_{b \in B} v_{bt} \boldsymbol z_{bts}\right) \label{1b_transportation_port_port} \\
\st \quad & \sum\nolimits_{s \in S} \boldsymbol z_{bts} = 1 && \hspace{-3.5cm} \forall b \in B, t \in T \label{1b_one_origin_port} \\
& \sum\nolimits_{h \in B} \boldsymbol y_{bsh} \le 1 && \hspace{-3.5cm} \forall b \in B, s \in S \label{1b_one_hub} \\
& \boldsymbol y_{bs0} + \sum\nolimits_{h \in B} \boldsymbol y_{bsh} \ge 1 && \hspace{-3.5cm} \forall b \in B, s \in S \label{1b_direct_or_hub} \\
& \boldsymbol y_{bsh} \leq \boldsymbol x_h && \hspace{-3.5cm} \forall b \in B, s \in S, h \in B \label{1b_selected_hubs} \\
& \boldsymbol y_{hs0} \geq \boldsymbol x_h && \hspace{-3.5cm} \forall h \in B, s \in S \label{1b_connections_from_hubs}
\end{align}

Constraint~\eqref{1b_one_origin_port} ensures that
exactly one origin port is used on each connection from a branch to a destination port.
In a similar manner,
constraint~\eqref{1b_one_hub} ensures that
at most one hub is selected on each connection from a branch to an origin port,
while constraint~\eqref{1b_direct_or_hub} makes sure that
at least one hub (\ie, exactly one, together with~\eqref{1b_one_hub}) is selected,
unless we use only direct transport on this connection anyway, \ie, $\boldsymbol y_{bs0} = 1$.
(Note that $\boldsymbol y_{bs0}$ is a continuous variable, while $\boldsymbol y_{bsh}$ is binary.)
By constraint~\eqref{1b_selected_hubs},
only branches that are selected as hubs may be used as hubs on a connection.
Finally,
constraint~\eqref{1b_connections_from_hubs} makes sure that connections starting in a hub use direct transport only.

The objective function consists of terms describing the set-up cost for hubs~\eqref{1b_set-up},
the consolidation cost at hubs~\eqref{1b_consolidation_hubs}
and at origin ports~\eqref{1b_consolidation_ports},
as well as the transportation cost from branches to origin ports~\eqref{1b_transportation_branch_port},
from branches to hubs~\eqref{1b_transportation_branch_hub},
and from origin ports to destination ports~\eqref{1b_transportation_port_port}.
For every origin port~$s \in S$,
the term $\sum_{b \in B} \sum_{t \in T} v_{bt} \boldsymbol z_{bts}$ in~\eqref{1b_consolidation_ports}
describes the volume that is shipped via $s$:
It is the sum of the volumes~$v_{bt}$ of all shipments that are routed via $s$, as determined by the variables $\boldsymbol z_{bts}$.
Similarly, for every hub~$h \in B$,
the term $\sum_{b \in B} \sum_{s \in S} \sum_{t \in T} v_{bt} \boldsymbol z_{bts} \boldsymbol y_{bsh} (1 - \boldsymbol y_{bs0})$ in~\eqref{1b_consolidation_hubs}
sums all volumes that are routed via the hub~$h$.
Note that the volume $v_{bt}$ of a shipment can be routed via $h$ if $h$ is selected as the hub for the route from branch~$b \in B$ to origin port~$s \in S$,
\ie, if $\boldsymbol y_{bsh} = 1$,
where $s$ is the origin port that is selected for the shipment via the variable $\boldsymbol z_{bts}$
(which is unique by~\eqref{1b_one_origin_port}).
This results in a total volume of $\sum_{b \in B} \sum_{s \in S} \sum_{t \in T} v_{bt} \boldsymbol z_{bts} \boldsymbol y_{bsh}$.
Finally, in each of the summands, only a percentage of $1 - \boldsymbol y_{bs0}$ of this volume is actually routed via the hub $h$
because the fraction $\boldsymbol y_{bs0}$ uses direct transport instead.
In a similar fashion,
the arguments of the cost functions in~\eqref{1b_transportation_branch_port}, \eqref{1b_transportation_branch_hub}, and \eqref{1b_transportation_port_port}
determine the volumes that are transported on the (direct) connections
from a branch $b \in B$ to an origin port $s \in S$,
from a branch $b \in B$ to a hub $h \in B$,
and from an origin port $s \in S$ to a destination port $t \in T$,
respectively.

%% file: chapters/3SchenkerSolution.tex
\section{The logistics provider's own solution} \label{chap_schenker_sol}

The logistics provider DB Schenker solved a simplified version of the problem in-house.
In terms of the mathematical problem, as presented in the previous section, the simplifications are:
\begin{itemize}
\item
The problem is split among the destination ports~$t\in T$:
Routing through the network and the decisions which branches to upgrade to hubs are made per destination port.
Parts of the consolidation effects (volume for different destination ports touching the same hinterland infrastructure) in land and ocean transport are thus ignored.
\item
An iterative heuristic based on an exhaustive search is used, instead of an exact method for the model.

\end{itemize}

\paragraph{Cost data}
As described in Section~\ref{chap_challenges}, the ocean transport costs are given per container.
In case there are different rates for the same relation, \eg, when both a CFS and an NVOCC are present, the lowest rate found in the rate database of the logistics provider is taken.
As land transport costs, the actual rates from the logistics provider, in this case DB Schenker's System Freight Network, are taken.
They are available up to 24 tons which constitute the full truck costs.

\paragraph{Shipment data}
A full year of historic shipment data is taken into account.
It is available on shipment level with shipment's origin branch, destination port, and the chargeable volume.
Since land transport costs are weight-dependent and shipment data is available in terms of chargeable volume, a dimensional factor of $300\kgpercbm$ is used to obtain the chargeable weight of a shipment.

To allow for result comparison, the shipment data also includes the actual origin ports.

\paragraph{Solution approach for the simplified problem}
For each destination port~$t\in T$, the simplified problem is solved by an iterative process switching between deciding which origin ports to use and which branches to upgrade to hubs:
Every iteration consists of two steps.
Let~$\bar b$ be a fixed number.
In the first step, the decision which origin ports to use is fixed and via an exhaustive search on the set~$\{H\subseteq B : |H|\leq\bar b\}$ of all possible sets of branches to upgrade to hubs, the optimal set of hubs is determined.
For a given set~$H$ of potential hubs, the costs themselves are evaluated exhaustively:
for each pair of origin-branch and origin-port all shipments are routed either directly or over one of the hubs~$h\in H$.
In the second step, the set of hubs is fixed and the origin ports to use are determined following the same logic as for the hubs.
The iterative process stops when no further cost reduction is achieved.

\paragraph{Solution quality and feasibility}
The split among destination ports~$t\in T$ results in a two-stage algorithm where a solution is computed per port and the results are simply merged:
A branch is upgraded to a hub when it is used as hub for at least one destination port.
The logistics provider's own solution
\begin{itemize}
\item
allows origin branches to use different hubs for transporting shipments to the same origin port, as different hubs can be selected for different destination ports~$t\in T$, which violates a constraint;
\item
is not able to guarantee that per port~$t\in T$ the optimal solution is found.
\end{itemize}
The logistics provider has thus decided to abandon its own solution and join forces with the co-authors of this paper to work on the solution approach presented in the next section.

%% file: chapters/4Hinterland.tex
\section{Combined hub and routing} \label{chap_tud_sol}

This section will show how the cost functions $\Gamma$ in the mixed integer program presented in Section~\ref{chap_math_model} can be modeled. First, land transport will be considered and then sea transport. In addition, the model together with a linearization is presented.

\subsection{Hinterland transport} \label{sec_modeling_hinterland_cost}

\begin{figure}[b]\centering
	\begin{minipage}[b]{.45\linewidth}
		\includegraphics[width=\linewidth]{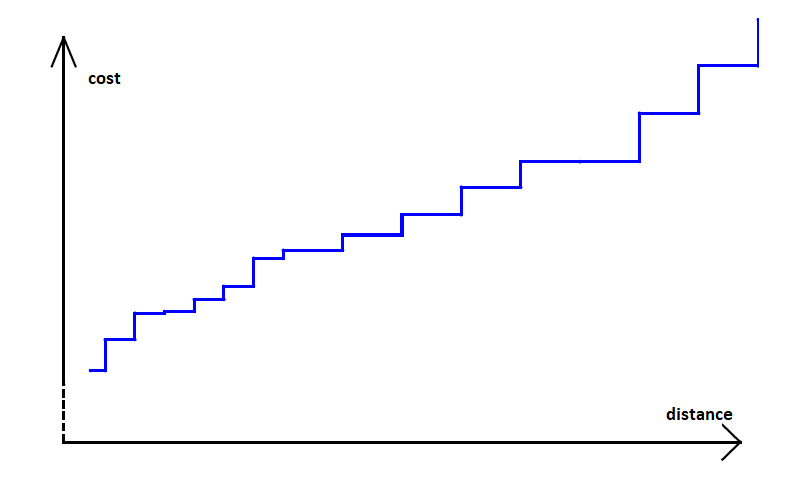}
		\caption{$\Gamma\Land$ with respect to distance, for a fixed volume}\label{fig_Gamma_land_distance}
	\end{minipage}
	\hspace{.08\linewidth}
	\begin{minipage}[b]{.45\linewidth}
		\includegraphics[width=\linewidth]{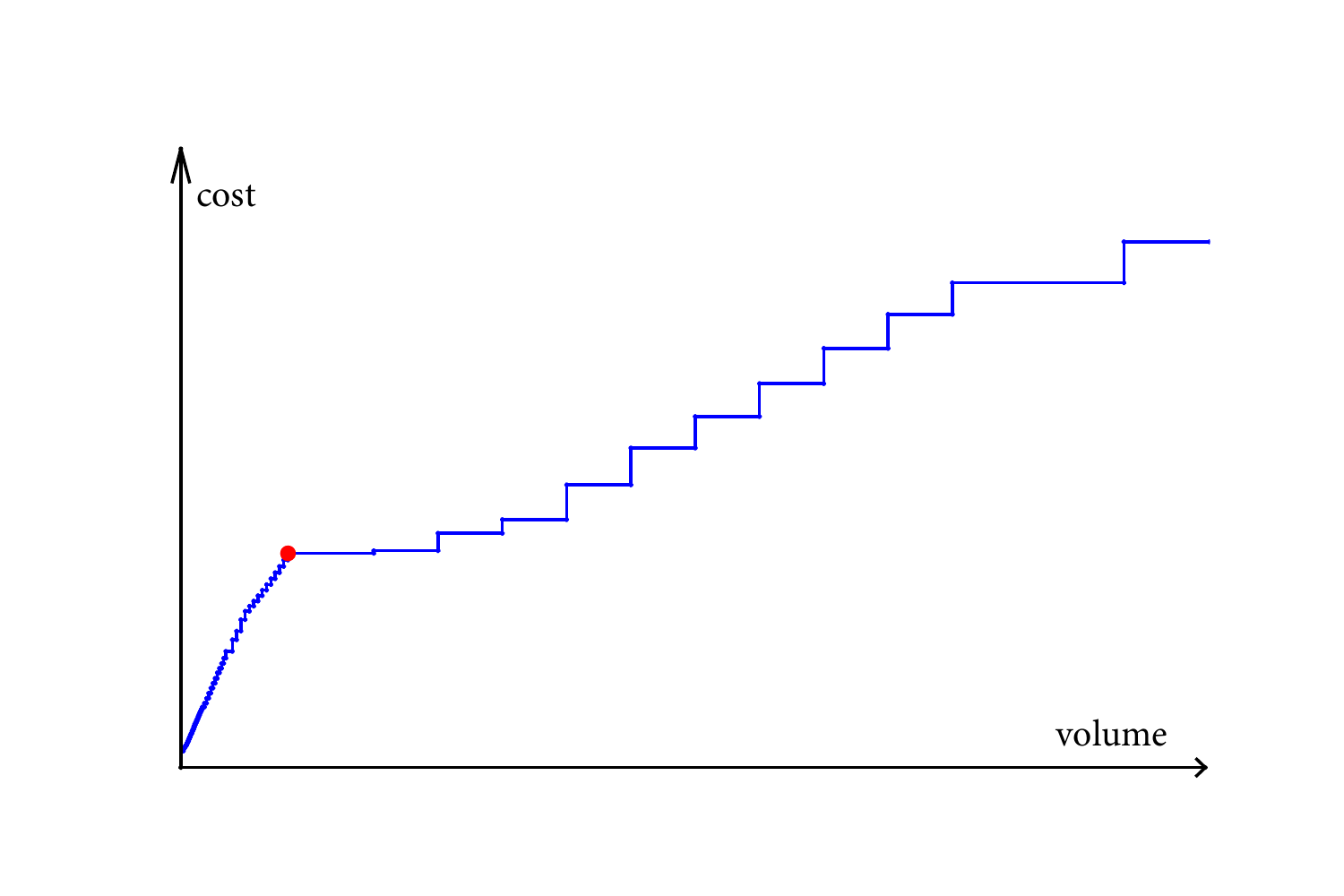}
		\caption{$\Gamma\Land$ with respect to volume, for a fixed distance} \label{fig_Gamma_land_volume}
	\end{minipage}
\end{figure}

DB Schenker has created a cost matrix $C$ for land transport based on their cost calculations. This matrix shows the cost for a transport depending on the shipment volume and distance. Let $u\Land_\cont$ denote the volume of a full container.
The transport costs can be expressed as a function $\Gamma\Land_{br}(v) = n \cdot C_{br}(u\Land_\cont) + C_{br}(u)$, where $v = n \cdot u\Land_\cont + u$ with $n$ maximal, \ie, $n$ is the maximal number of containers that can be filled completely and $u$ is the rest volume that does not fit into these containers. Here, $C_{br}(x)$ denotes the value in the matrix corresponding to the distance between a branch $b \in B$ and a branch or origin port $r \in B \cup S$, as well as volume~$x$. The function $\Gamma\Land_{br}(v)$ is non-decreasing and piecewise constant. It jumps whenever $u$ changes from one volume range to the next or when the number~$n$ of full containers changes. Figure~\ref{fig_Gamma_land_distance} shows an example for the relation between the distance and the cost for a fixed volume. In Figure~\ref{fig_Gamma_land_volume}, the distance is fixed and the costs are given with respect to the volume.

The cost functions have $63$ different volume ranges for different price structures. To use the hinterland transport cost functions $\Gamma\Land_{br}(v)$ in the proposed mathematical programs, for every fixed distance, the first constant pieces were approximated by a linear function for $u<u\Land_\cont/10$. The reason for this is that the formulation was too complex to solve otherwise and the cost functions do not differ much from a linear shape in the beginning; see Figure~\ref{fig_Gamma_land_volume}, where the red point is at the volume $u\Land_\cont/10$. For higher values, an exact modeling approach was used. Denote the corresponding approximation of $\Gamma\Land_{br}(v)$ by $\widetilde{\Gamma}\Land_{br}(v)$ and the $v$-values of the breakpoints or jumps of $\widetilde{\Gamma}\Land_{br}(v)$ by $v^{(0)} = u\Land_\cont/10,v^{(1)},\ldots,v^{(j)} = u\Land_\cont$.

Let $\widetilde{C}^{(i)}_{br} = \widetilde{\Gamma}\Land_{br}(v^{(i)})$ for $i \in \{0, \ldots, j\}$ be the corresponding function values. We assume $\widetilde{C}^{(i)}_{br}$ to be the constant value in the range $(v^{(i - 1)}, v^{(i)}]$ for $i \in \{1, \ldots, j\}$. For the range $[0, v^{(0)}]$, the value of $\widetilde{\Gamma}\Land_{br}(v)$ is modeled as a linear function between the points $(0,0)$ and $(v^{(0)}, \widetilde{C}^{(0)}_{br})$. The same pattern repeats for the second container and so on.

The approximation $\widetilde{\Gamma}\Land_{br}(v)$ can be modeled in the mathematical program by introducing the following new variables:
\begin{itemize}
	\item $\boldsymbol n_{br} \in \mathbb{Z}_{\geq 0}$ for all $b \in B$ and $r \in B \cup S$, determining the number of full containers between branch $b \in B$ and branch or origin port $r \in B \cup S$,
	\item $\boldsymbol u^{(0)}_{br} \in [0,1]$ for all $b \in B$ and $r \in B \cup S$, determining the volume not fitting into the full containers between a branch $b \in B$ and a branch or origin port $r \in B \cup S$, expressed as a percentage of $v^{(0)}$, if this volume is at most $v^{(0)}$, otherwise $\boldsymbol u^{(0)}_{br}=0$,
	\item $\boldsymbol u^{(i)}_{br} \in \{0, 1\}$ for all $b \in B$, $r \in B \cup S$, and $i \in \{1, \ldots, j-1\}$, determining if the volume not fitting into the full containers between a branch $b \in B$ and a branch or origin port $r \in B \cup S$ lies in the range $(v^{(i - 1)}, v^{(i)}]$ ($\boldsymbol u^{(i)}_{br}=1$) or not ($\boldsymbol u^{(i)}_{br}=0$).
\end{itemize}
Note that there is no variable $\boldsymbol u^{(j)}_{br}$ since the last step in the cost functions also corresponds to the transport costs of a full container, which is modeled through the variable $\boldsymbol n_{br}$.

The approximate hinterland transport costs for a branch $b \in B$ and a branch or origin port $r \in B \cup S$ are then determined by the expression
\begin{equation*}
	\widetilde{C}^{(j)}_{br} \boldsymbol n_{br} + \sum\nolimits_{i = 0}^{j-1} \widetilde{C}^{(i)}_{br} \boldsymbol u^{(i)}_{br}\,,
\end{equation*}
which can be used to replace the function $\Gamma\Land_{br}(v)$ in the terms~\eqref{1b_transportation_branch_port} and~\eqref{1b_transportation_branch_hub}.

Moreover, in order to obtain the desired values of the new variables, we add the following constraints to the model from Section~\ref{chap_math_model}:
\begin{align*}
&\begin{array}{l}
\displaystyle \sum\nolimits_{t \in T} \left(v_{bt} \boldsymbol z_{bts} \boldsymbol y_{bs0} + \sum\nolimits_{c \in B} v_{ct} \boldsymbol z_{cts} \boldsymbol y_{csb}(1-\boldsymbol y_{cs0})\right)\\
\displaystyle \leq v^{(j)} \boldsymbol n_{bs} + \sum\nolimits_{i = 0}^{j-1} v^{(i)} \boldsymbol u^{(i)}_{bs}
\end{array}
&& \forall b \in B, s \in S \\
	& \sum_{s \in S} \sum_{t \in T} v_{bt} \boldsymbol z_{bts} \boldsymbol y_{bsh}(1-\boldsymbol y_{bs0}) \leq v^{(j)} \boldsymbol n_{bh} + \sum_{i = 0}^{j-1} v^{(i)} \boldsymbol u^{(i)}_{bh} && \forall b, h \in B
\end{align*}

In order to ensure that at most one of the variables $\boldsymbol u^{(i)}_{br}$ takes a positive value for each $b \in B$ and $r \in B \cup S$, we finally add the following constraint:
\begin{align*}
	& \sum\nolimits_{i = 0}^{j-1} \boldsymbol u^{(i)}_{br} \leq 1 & \forall b \in B, r \in B \cup S
\end{align*}

An exact model can be obtained by using a binary variable $\boldsymbol u^{(i)}_{br}$ for every jump of the original function and no continuous variable $\boldsymbol u^{(0)}_{br}$.

\subsection{Sea transport} 
Any transport below $40\cbm$ can be transported either using own containers (FCL) or by an NVOCC if both services are available. A volume between $40$ and $55\cbm$ is always transported by FCL. The volume of $55\cbm$ is a full container such that any volume above $55\cbm$ is split into full containers, and the remaining volume is treated respectively. 

DB Schenker has provided cost calculations for NVOCC and FCL transport for every possible connection from an origin port $s \in S$ to a destination port $t \in T$.
The sea transport cost function can be written as $\Gamma\Sea_{st}(v) = n \cdot C\Sea_{st\FCL} + u \cdot C\Sea_{st\NVOCC}$, where $n \cdot u\Sea_\cont + u \geq v$ and $u \leq u^{\lim}_{st}$, \ie, $n$ is the required number of full containers (of volume $u\Sea_\cont = 55\cbm$) and $u$ is the rest volume that does not fit into these containers (if this is at most $u \leq u^{\lim}_{st}$). The numbers $C\Sea_{st\FCL}$ and $C\Sea_{st\NVOCC}$ describe the FCL (per container) and NVOCC (per m$^3$) transportation costs between $s$ and $t$ in the provided tables, respectively. The value of  $u^{\lim}_{st}$ corresponds to the minimum of $u_{\max} = 40\cbm$ and the NVOCC volume needed to match the costs of an FCL shipment, which can be written as $u^{\lim}_{st} = \min \left\lbrace u_{\max}, \frac{C\Sea_{st\FCL}}{C\Sea_{st\NVOCC} }\right\rbrace$. However, using $u^{\lim}_{st} = u_{\max}$ instead leads to the same results due to the minimization. Figure~\ref{fig_Gamma_sea} illustrates this cost function.

\begin{figure}[t]
	\centering
	\begin{tikzpicture}[xscale=0.03,yscale=0.03]
		\draw[->] (0,0) -- (120,0) node[anchor=west] {volume [m$^3$]};
		\draw (22.4,0) node {$\scriptstyle{+}$} node[below] {$u_{st}^{\lim}$};
		\draw (55,0) node {$\scriptstyle{+}$} node[below] {$u\Sea_\cont$};
	
		\draw[->] (0,0) -- (0,60) node[anchor=south] {cost};

		\draw (0,0) -- (22.4,26.282) -- (55,26.282) -- (77.4,52.564) -- (110,52.564);
	\end{tikzpicture}
\caption{$\Gamma\Sea$ with respect to volume} \label{fig_Gamma_sea}
\end{figure}

The sea transport cost functions $\Gamma\Sea_{st}(v)$ can be modeled exactly in the mathematical program using the following new variables:
\begin{itemize}
	\item $\boldsymbol n_{st} \in \mathbb{Z}_{\geq 0}$ for all $s \in S$ and $t \in T$, determining the number of full containers between origin port $s \in S$ and destination port $t \in T$,
	\item $\boldsymbol u_{st} \in [0, u^{\lim}_{st}]$ for all $s \in S$ and $t \in T$, determining the amount shipped by an NVOCC (in m$^3$) between origin port $s \in S$ and destination port $t \in T$.
\end{itemize}

The term \eqref{1b_transportation_port_port} in the objective function is replaced by
\begin{equation*}
	\sum\nolimits_{s \in S} \sum\nolimits_{t \in T} \left(C\Sea_{st\FCL} \boldsymbol n_{st} + C\Sea_{st\NVOCC} \boldsymbol u_{st}\right).
\end{equation*}
Moreover, the following constraint is added, which ensures the correct distribution of the volume among the new variables:
\begin{align*}
	& \sum\nolimits_{b \in B} v_{bt} \boldsymbol z_{bts} \leq u\Sea_\cont \boldsymbol n_{st} + \boldsymbol u_{st} & \forall s \in S, t \in T
\end{align*}

In case there is only an FCL cost value given for a connection between an origin port $s \in S$ and a destination port $t \in T$, the corresponding variable $\boldsymbol u_{st}$ is omitted (or, equivalently, set to 0) in the model and the new constraints.

Suppose there is only an NVOCC cost value given. Then the variable $\boldsymbol n_{st}$ can be omitted, leading to a model where the bound $u_{\max}$ is still strictly enforced. This might make the model infeasible, which can be avoided by allowing larger volumes in this case, but only while paying a high penalty for them. This can be modeled by not omitting the variable $\boldsymbol n_{st}$, but leaving the model as it is and just setting the coefficient in front of $\boldsymbol n_{st}$ in the objective function to a high penalty value (depending on the other data, \eg, $10^8$) instead of $C\Sea_{st\FCL}$.

\subsection{Summary of the model}

In the following, we present the complete model description, including the modifications regarding the cost functions described in this section.
\begin{itemize}
\item $\boldsymbol z_{bts} \in \{0, 1\}$ for all $b \in B$, $t \in T$, and $s \in S$, determining if $s$ is used as origin port (1) or not (0) on the connection from $b$ to $t$,
\item $\boldsymbol x_b \in \{0, 1\}$ for all $b \in B$, determining if $b$ is a hub ($1$) or not ($0$),
\item $\boldsymbol y_{bs0} \in [0, 1]$ for all $b \in B$ and $s \in S$, determining which percentage of the volume on the connection from $b$ to $s$ uses direct transport,
\item $\boldsymbol y_{bsh} \in \{0, 1\}$ for all $b \in B$, $s \in S$, and $h \in B$, determining if $h$ is used as a hub ($1$) or not ($0$) on the connection from $b$ to $s$,
\item $\boldsymbol n_{br} \in \mathbb{Z}_{\geq 0}$ for all $b \in B$ and $r \in B \cup S$, determining the number of full containers between branch $b \in B$ and branch or origin port $r \in B \cup S$,
\item $\boldsymbol u^{(0)}_{br} \in [0,1]$ for all $b \in B$ and $r \in B \cup S$, determining the volume not fitting into the full containers between a branch $b \in B$ and a branch or origin port $r \in B \cup S$, expressed as a percentage of $v^{(0)}$, if this volume is at most $v^{(0)}$, otherwise $\boldsymbol u^{(0)}_{br}=0$,
\item $\boldsymbol u^{(i)}_{br} \in \{0, 1\}$ for all $b \in B$, $r \in B \cup S$, and $i \in \{1, \ldots, j-1\}$, determining if the volume not fitting into the full containers between a branch $b \in B$ and a branch or origin port $r \in B \cup S$ lies in the range $(v^{(i - 1)}, v^{(i)}]$ ($\boldsymbol u^{(i)}_{br}=1$) or not ($\boldsymbol u^{(i)}_{br}=0$),
\item $\boldsymbol n_{st} \in \mathbb{Z}_{\geq 0}$ for all $s \in S$ and $t \in T$, determining the number of full containers between origin port $s \in S$ and destination port $t \in T$,
\item $\boldsymbol u_{st} \in [0, u^{\lim}_{st}]$ for all $s \in S$ and $t \in T$, determining the amount shipped by an NVOCC (in m$^3$) between origin port $s \in S$ and destination port $t \in T$.
\end{itemize}
\begin{align}
	\min \quad & \sum_{h \in B} e_h \boldsymbol x_h + \sum_{s \in S} g_s \sum_{b \in B} \sum_{t \in T} v_{bt} \boldsymbol z_{bts} \notag \\
	& + \sum_{h \in B} f_{h} \sum_{b \in B} \sum_{s \in S} \sum_{t \in T} v_{bt} \boldsymbol z_{bts} \boldsymbol y_{bsh}(1-\boldsymbol y_{bs0}) \notag \\
	& + \sum_{b \in B} \sum_{s \in S} (\widetilde{C}^{(j)}_{bs} \boldsymbol n_{bs} + \sum_{i = 0}^{j-1} \widetilde{C}^{(i)}_{bs} \boldsymbol u^{(i)}_{bs}) \notag \\
	& + \sum_{b \in B} \sum_{h \in B} (\widetilde{C}^{(j)}_{bh} \boldsymbol n_{bh} + \sum_{i = 0}^{j-1} \widetilde{C}^{(i)}_{bh} \boldsymbol u^{(i)}_{bh}) \notag \\
	& + \sum_{s \in S} \sum_{t \in T} \left(C\Sea_{st\FCL} \boldsymbol n_{st} + C\Sea_{st\NVOCC} \boldsymbol u_{st}\right) \notag \\
  \st  \quad & \sum\nolimits_{s \in S} \boldsymbol z_{bts} = 1 && \forall b \in B, t \in T \notag\\
  	& \sum\nolimits_{h \in B} \boldsymbol y_{bsh} \leq 1 && \forall b \in B, s \in S \notag\\
	& \boldsymbol y_{bs0} + \sum\nolimits_{h \in B} \boldsymbol y_{bsh} \geq 1 && \forall b \in B, s \in S \label{weg1}\\
	& \boldsymbol y_{bsh} \leq \boldsymbol x_h && \forall b, h \in B, s \in S\notag\\
	& \boldsymbol y_{hs0} \geq \boldsymbol x_h && \forall h \in B, s \in S \label{weg2}\\
	& \begin{array}{l}\displaystyle\sum_{t \in T} \left(v_{bt} \boldsymbol z_{bts} \boldsymbol y_{bs0} + \sum_{c \in B} v_{ct} \boldsymbol z_{cts} \boldsymbol y_{csb}(1-\boldsymbol y_{cs0})\right)\\
	\displaystyle\leq v^{(j)} \boldsymbol n_{bs} + \sum\nolimits_{i = 0}^{j-1} v^{(i)} \boldsymbol u^{(i)}_{bs}\end{array} && \forall b \in B, s \in S \notag\\
	& \sum_{s \in S} \sum_{t \in T} v_{bt} \boldsymbol z_{bts} \boldsymbol y_{bsh} (1{-}\boldsymbol y_{bs0})\leq v^{(j)} \boldsymbol n_{bh} + \!\!\sum_{i = 0}^{j-1} v^{(i)} \boldsymbol u^{(i)}_{bh} && \forall b, h \in B \notag\\
	& \sum\nolimits_{i = 0}^{j-1} \boldsymbol u^{(i)}_{br} \leq 1 && \forall b \in B, r \in B \cup S \notag\\
	& \sum\nolimits_{b \in B} v_{bt} \boldsymbol z_{bts} \leq u\Sea_\cont \boldsymbol n_{st} + \boldsymbol u_{st} && \forall s \in S, t \in T\notag
\end{align}
This model is now a mixed-integer cubic program.
We propose a linearized version of the cubic terms
by introducing the following new variables:
\begin{itemize}
\item $\boldsymbol v_{bs0} \in \mathbb{R}_{\geq 0}$ for all $b \in B$ and $s \in S$, determining the volume that is sent from branch $b$ to origin port $s$ using direct transport,
\item $\boldsymbol v_{bsh} \in \mathbb{R}_{\geq 0}$ for all $b \in B$, $s \in S$, and $h \in B$, determining the volume that is sent from branch $b$ to origin port $s$ via the hub $h$.
\end{itemize}
The idea is that $\boldsymbol v_{bsh} = \sum_{t \in T} v_{bt} \boldsymbol z_{bts} \boldsymbol y_{bsh}(1-\boldsymbol y_{bs0})$ for $b \in B$, $s \in S$, and $h \in B$, and $\boldsymbol v_{bs0} = \sum_{t \in T} v_{bt} \boldsymbol z_{bts} \boldsymbol y_{bs0}$ for $b \in B$ and $s \in S$ are supposed to hold. In order to ensure this, we add the following constraints to the model:
\begin{align}
  & \sum\nolimits_{t \in T} v_{bt} \boldsymbol z_{bts} = \boldsymbol v_{bs0} + \sum\nolimits_{h \in B} \boldsymbol v_{bsh} && \forall b \in B, s \in S \label{1b_split_volume} \\
  & \boldsymbol v_{bsh} \leq M \boldsymbol y_{bsh} && \forall b \in B, s \in S, h \in B, \notag
\end{align}
with $M$ being an upper bound of the value of $\boldsymbol v_{bsh}$, \eg, $M_b = \sum_{t \in T} v_{bt}$ for every $b \in B$. The correctness then follows from~\eqref{1b_one_hub}.
Replacing the ``$=$'' in~\eqref{1b_split_volume} by ``$\leq$'' leads to an equivalent model due to the optimization.
In addition, the variable $\boldsymbol y_{bs0}$ can be eliminated everywhere and the constraint \eqref{weg1} can be neglected. The constraint \eqref{weg2} can then be replaced by 
\begin{align*}
	&\boldsymbol y_{hsc}\leq 1-\boldsymbol x_h &  \forall h \in B, s \in S, c \in B .
\end{align*}
We can now replace all cubic terms in the objective function and the constraints by the new variables,
which results in a mixed-integer linear program.

Note that also other linearizations of the cubic terms are possible (and maybe worth investigating further), but this one has the advantage that it uses only $O(|B|^2\cdot|S|)$ many new variables, in contrast to the standard linearization of the products $\boldsymbol z_{bts} \boldsymbol y_{bsh}$, which would need $O(|B|^2\cdot|S|\cdot|T|)$ many new variables. This is possible by exploiting the constraint~\eqref{1b_one_hub}, similarly to known linearizations of the standard hub location problem; see, \eg, \cite{KaraTansel}.

%% file: chapters/5Results.tex
\section{Results}\label{chap_results}%

The hinterland routing of the trucks depends on the ocean routing of the containers. 
This fact cannot be considered by the two-stage approach established in Section~\ref{chap_schenker_sol}.
In Section~\ref{chap_tud_sol}, a mathematical model for the non-partitioned solution approach is presented.
It yields significantly better results, yet it is coupled with substantially more computational effort.

For implementing and solving the mathematical model (in the linearized version) described in Sections~\ref{chap_math_model} and~\ref{chap_tud_sol}, we used \CC{} together with the optimization software \cite{gurobi}.
The transportation costs were modeled as presented in Section~\ref{chap_tud_sol}, and the set-up and consolidation costs were omitted.
From the historic shipment data, as described in Section~\ref{chap_schenker_sol}, two classes of seven instances each were generated.
The first class corresponds to shipments originating in Spain and Portugal.
Each instance consists of the shipments sent in a specific week.
For each of them, the number of different relations (origin branch--destination port) on which shipments are to be routed lies between 35 and 75.
The second class of instances corresponds to shipments originating in Spain, Portugal, and France.
Here, each instance covers between 100 and 200 different relations.

In order to assess the quality of the exact modeling approach, we compare the results to the optimal solutions of the restricted network planning problem in which no hubs are built for consolidating cargo in the hinterland, but instead direct transport is always used between the origin branch of a shipment and its origin port.
Hence, the task of this problem is only to select an origin port for every shipment.
The resulting model is significantly smaller than the one with hubs and all instances could be solved optimally in at most one minute.

Regarding the model with hubs, within a time limit of one hour for each instance, Gurobi could solve one instance optimally, the other instances of the first class up to a gap of between 1.22\% and 3.83\%, and the ones of the second class up to a gap of between 3.28\% and 12.53\%.
The solutions of our model send between 20.79\% and 49.65\%, with an average of 38.29\%, of the total shipment volume via hubs, while the remaining volume is sent directly from its origin branch to its origin port.
In all instances except for one of the second class, the best solution obtained with our model within one hour has a lower total cost than the optimal solution without hubs.
The improvement for these instances lies between 3.24\% and 7.11\% for the first class, and between 1.13\% and 2.94\% for the second class.
Sometimes, the land transport costs increase, but the sea transport costs always decrease significantly.
This holds both when comparing the costs the model uses, including the linear approximation of parts of the land transport cost functions, and when comparing the exact costs instead.

In conclusion, the routing decisions in the hub location problem should be made at once. This also includes the decision whether a hub is opened or not if this is part of the problem description. Mathematical models can represent the problem statement and can be solved.

It would be desirable to decrease the computational effort for solving the mathematical model in order to achieve similar results also for larger instances.
For example, all shipments originating in Europe in a given week result in an instance covering around 1000 different relations.
Future work could establish a Lagrangian decomposition method for the model based on \cite{RostamiBuchheim}.
Moreover, an even more realistic model could be achieved by taking into account that the shipment volumes are usually not known precisely in advance for a longer time horizon. This uncertainty can be modeled using a stochastic optimization approach.